\documentclass[eqno]{amsart}
\usepackage[dvips,final]{epsfig}
\usepackage{amsmath}
\usepackage{amsfonts}
\usepackage{amssymb}
\usepackage{latexsym}
\usepackage{color}
\title[Whitney umbrellas and swallowtails]
{Whitney umbrellas 
and 
swallowtails}
\dedicatory{
Dedicated to Professor Shyuichi Izumiya on the occasion of his sixtieth birthday 
}
\author[
T. Nishimura
]{
Takashi Nishimura
}
\address{
Research Group of Mathematical Sciences,  
Research Institute of Environment and Information Sciences, 
Yokohama National University, 
Yokohama 240-8501, JAPAN. 
}
\email{nishimura-takashi-yx@ynu.jp}

\begin{document}
\begin{abstract}\quad 
In this paper, we introduce the notions of map-germs of pedal unfolding type 
and normalized Legendrian map-germs; and then we show that the fundamental theorem of calculus provides 
a natural one to one correspondence between 
Whitney umbrellas of pedal unfolding type and normalized swallowtails.   
\end{abstract}
\date{}
\maketitle
\noindent

\newtheorem{theorem}{Theorem} 
\newtheorem{corollary}{Corollary} 
\newtheorem{lemma}{Lemma}[section]
\newtheorem{proposition}{Proposition}
\newtheorem{definition}{Definition}
\newtheorem{example}{Example}[section] 
\newtheorem{conjecture}{Conjecture}
\def\theenumi{\roman{enumi}}
\noindent
\section{Introduction} 
\label{intro}
\noindent 
The map-germ 
\begin{equation}
f(x,y)=(xy, x^2, y) 
\label{Whitney}   
\end{equation}
is well-known as the normal form of Whitney umbrella after Whitney's pioneer works \cite{whitney, whitney2}.    
For the map-germ~\eqref{Whitney}, compose the following two coordinate transformations:   
$h_s(x,y)=(x, x^2+y)$ and $h_t(X,Y,Z)=(X,-Z,-Y+Z)$ where 
$(X,Y,Z)$ is the standard coordinates of the target space $\mathbb{R}^3$.     
Then, we have the following:   
\begin{equation} \label{g}
g(x,y)=h_t\circ f\circ h_s(x,y)=(x^3+xy, -x^2-y, y). 
\end{equation}
Put 
{\small 
\begin{equation} \label{G}
{G}(x,y)  =  
\left(\int_0^x (x^3+xy)dx, \int_0^x (-x^2-y)dx, y\right) 
 =  \left(\frac{1}{4}x^4+\frac{1}{2}x^2y, -\frac{1}{3}x^3-xy, y\right).   
\end{equation}
}
For the map-germ~\eqref{G}, compose the following two scaling transformations:   
${H}_s(x,y)=(x,\frac{1}{6}y)$ and ${H}_t(X,Y,Z)=(12X,12Y,6Z)$.     
Then, we have the following:   
\begin{equation} \label{swallowtail}
{H}_t\circ {G}\circ {H}_s(x,y)=(3x^4+x^2y,-4x^3- 2xy, y).
\end{equation}
The map-germ~\eqref{swallowtail} is well-known as the normal form of swallowtail (for instance, see 
\cite{brucegiblin} p.129).   
\par
Two $C^\infty$ map-germs $\varphi,\psi:(\mathbb{R}^2,0)\to (\mathbb{R}^3,0)$ are said to be 
$\mathcal{A}$-equivalent if there exist germs of $C^\infty$ diffeomorphisms 
$h_s: (\mathbb{R}^2,0)\to (\mathbb{R}^2,0)$ and  $h_t: (\mathbb{R}^3,0)\to (\mathbb{R}^3,0)$ 
such that $\psi=h_t\circ \varphi\circ h_s$.       
A $C^\infty$ map-germ $\varphi:(\mathbb{R}^2,0)\to (\mathbb{R}^3,0)$ is called 
a {\it Whitney umbrella} (resp., {\it swallowtail}) if $\varphi$ is $\mathcal{A}$-equivalent to 
\eqref{Whitney} (resp., \eqref{swallowtail}).       
As above, the Whitney umbrella~\eqref{Whitney} produces the swallowtail~\eqref{swallowtail} via \eqref{g} and \eqref{G}.     
By the converse procedure, the swallowtail~\eqref{swallowtail} produces the Whitney umbrella~\eqref{Whitney}.    
\par 
Note that it is impossible to produce a swallowtail by integrating \eqref{Whitney} directly.       
This is because the discriminant set of \eqref{swallowtail} is not diffeomorphic to the discriminant set of the 
following \eqref{5}:       
\begin{equation} \label{5}
(x,y)\mapsto 
\left(
\int_0^x xy dx, \int_0^x x^2dx, y 
\right). 
\end{equation}
Note further that 
the form \eqref{g} may be written as follows:   
\[
g(x,y)=\left(x(x^2+y), -(x^2+y), y\right)
=\left(b \left( -x, -(x^2+y)\right), y\right), 
\]
where $b(X,Y)=(XY,Y)$ ($b$ stands for \lq\lq the blow down\rq\rq ).
\begin{definition} \label{definition 1}
{\rm     
\begin{enumerate}
\item 
A $C^\infty$ map-germ $\varphi: (\mathbb{R}^2,0)\to (\mathbb{R}^3,0)$ 
having the following form \eqref{pedal} 
is said to be {\it of pedal unfolding type}.     
\begin{equation} \label{pedal}
\varphi(x,y)=\left(n(x,y)p(x,y), p(x,y), y\right)
=\left(b \left( n(x,y), p(x,y)\right), y\right) 
\end{equation}
where $n: (\mathbb{R}^2,0)\to (\mathbb{R},0)$ is a $C^\infty$ function-germ 
such that $\frac{\partial n}{\partial x}(0,0)\ne 0$ and 
$p: (\mathbb{R}^2,0)\to (\mathbb{R},0)$ is a $C^\infty$ function-germ.    
\item 
For a $C^\infty$ map-germ of pedal unfolding type 
$\varphi(x,y)=(n(x,y)p(x,y), p(x,y),y)$, put 
\[
\mathcal{I}(\varphi)=\left(
\int_0^x n(x,y)p(x,y)dx, \int_0^x p(x,y)dx, y
\right). 
\]
The map-germ $\mathcal{I}(\varphi): (\mathbb{R}^2,0)\to (\mathbb{R}^3,0)$ is 
called the {\it integration of }$\varphi$.    
\item 
A $C^\infty$ map-germ $\Phi: (\mathbb{R}^m,0)\to (\mathbb{R}^{m+1},0)$ is called 
a {\it Legendrian map-germ} if there exists a germ of $C^\infty$ vector field 
$\nu_\Phi : (\mathbb{R}^m,0)\to T_1\mathbb{R}^{m+1}$ along $\Phi$ 
such that the following conditions hold 
where the dot in the center stands for the scalar product of two vectors of $T_{\Phi(x,y)}\mathbb{R}^{m+1}$ and 
$T_1\mathbb{R}^{m+1}$ stands for the unit tangent bundle of $\mathbb{R}^{m+1}$.  
\begin{enumerate}
\item 
$\frac{\partial \Phi}{\partial x_1}(x_1, \ldots, x_m)\cdot \nu_\Phi(x_1, \ldots, x_m)=\cdots =$ 
$\frac{\partial \Phi}{\partial x_m}(x_1, \ldots, x_m)\cdot \nu_\Phi(x_1, \ldots, x_m)$ 
$=0$.  
\item The map-germ $L_\Phi: (\mathbb{R}^{m},0)\to T_1\mathbb{R}^{m+1}$ defined by 
\[
L_\Phi(x_1, \ldots, x_m)=(\Phi(x_1, \ldots, x_m), \nu_\Phi(x_1, \ldots, x_m))
\] 
is non-singular.      
The map-germ $L_\Phi$ is called a {\it Legendrian lift of} $\Phi$.     
\end{enumerate} 
The $C^\infty$ vector field $\nu_\Phi$ is called a {\it unit normal vector field of} $\Phi$.     
\item 
A Legendrian map-germ $\Phi:  (\mathbb{R}^2,0)\to (\mathbb{R}^3,0)$ 
is said to be {\it normalized} if $\Phi$ satisfies the following three conditions: 
\begin{enumerate}
\item $\Phi$ has the form $\Phi(x,y)=(\Phi_1(x,y), \Phi_2(x,y),y)$.  
\item $\frac{\partial \Phi_2}{\partial x}(0,0)=0$.   
\item $\nu_\Phi(0,0)$ is $\frac{\partial}{\partial X}$ or $-\frac{\partial}{\partial X}$.     
\end{enumerate}    
\item 
For a normalized Legendrian map-germ  
$\Phi(x,y)=(\Phi_1(x,y), \Phi_2(x,y),y)$, put 
\[
\mathcal{D}(\Phi)=\left(
\frac{\partial \Phi_1}{\partial x}(x,y), 
\frac{\partial \Phi_2}{\partial x}(x,y), y
\right). 
\]
The map-germ $\mathcal{D}(\Phi): (\mathbb{R}^2,0)\to (\mathbb{R}^3,0)$ is 
called the {\it differentiation of }$\Phi$.    
\end{enumerate}
}
\end{definition}
Since any map-germ $\varphi: (\mathbb{R}^2,0)\to (\mathbb{R}^3,0)$ which is 
a germ of one-parameter pedal unfolding of a spherical pedal curve 
has the form \eqref{pedal} (see \cite{nishimura3}), it is reasonable that  
a map-germ $\varphi$ having the form \eqref{pedal} is said to be of pedal unfolding type.     
As shown in \cite{nishimura3},  not only non-singular map-germs 
but also Whitney umbrellas may be realized 
as germs of one-parameter pedal unfoldings of spherical pedal curves.      
For details on Legendrian map-germs, see \cite{arnoldetall, izumiya,zakalyukin2, zakalyukin}.        
Note that both \eqref{G} and \eqref{cuspidal} are normalized Legendrian map-germs.        
\begin{proposition} \label{proposition 1}
\begin{enumerate}
\item For a $C^\infty$ map-germ of pedal unfolding type 
$\varphi : (\mathbb{R}^2,0)\to (\mathbb{R}^3,0)$, 
$\mathcal{I}(\varphi)$ is a normalized Legendrian map-germ.  
\item For a normalized Legendrian map-germ  
$\Phi : (\mathbb{R}^2,0)\to (\mathbb{R}^3,0)$, 
$\mathcal{D}(\Phi)$ is a map-germ of pedal unfolding type.  
\end{enumerate}
\end{proposition}
Put 
\begin{eqnarray*}
\mathcal{W} & = & 
\left\{ \varphi: (\mathbb{R}^2,0)\to (\mathbb{R}^3,0)
\; \mbox{ Whitney umbrella of pedal unfolding type}\right\}, \\ 
\mathcal{S} & = & 
\left\{ \Phi: (\mathbb{R}^2,0)\to (\mathbb{R}^3,0)\; \mbox{ normalized swallowtail}
\right\}.  
\end{eqnarray*}
\par 
The main purpose of this paper is to show the following Theorem~\ref{theorem 1}:   
\begin{theorem} \label{theorem 1}
\begin{enumerate}
\item The map $\mathcal{I}: \mathcal{W}\to \mathcal{S}$ 
defined by $\mathcal{W}\ni\varphi\mapsto \mathcal{I}(\varphi)\in\mathcal{S}$ 
is well-defined and bijective. 
\item The map $\mathcal{D}: \mathcal{S}\to \mathcal{W}$ 
defined by $\mathcal{S}\ni\Phi\mapsto \mathcal{D}(\Phi)\in\mathcal{W}$ 
is well-defined and bijective. 
\end{enumerate}     
\end{theorem}
Incidentally, we show the following Theorem~\ref{theorem 2}.    
A $C^\infty$ map-germ $\Phi:  (\mathbb{R}^2,0)\to (\mathbb{R}^3,0)$ is called a {\it cuspidal edge} if 
$\Phi$ is $\mathcal{A}$-equivalent to the following \eqref{cuspidal} : 
\begin{equation} \label{cuspidal}
(x,y)\mapsto \left(\frac{1}{3}x^3, \frac{1}{2}x^2, y\right). 
\end{equation}      
Put 
\begin{eqnarray*}
\mathcal{N} & = & 
\left\{ \varphi: (\mathbb{R}^2,0)\to (\mathbb{R}^3,0)\; \mbox{ non-singular map-germ of pedal unfolding type}\right\}, \\ 
\mathcal{C} & = & 
\left\{ \Phi: (\mathbb{R}^2,0)\to (\mathbb{R}^3,0)\; \mbox{ normalized cuspidal edge}
\right\}.  
\end{eqnarray*}
\begin{theorem} \label{theorem 2}
\begin{enumerate}
\item The map $\mathcal{I}: \mathcal{N}\to \mathcal{C}$ 
defined by $\mathcal{N}\ni\varphi\mapsto \mathcal{I}(\varphi)\in\mathcal{C}$ 
is well-defined and bijective.    
\item The map $\mathcal{D}: \mathcal{C}\to \mathcal{N}$ 
defined by $\mathcal{C}\ni\Phi\mapsto \mathcal{D}(\Phi)\in\mathcal{N}$ 
is well-defined and bijective.    
\end{enumerate}     
\end{theorem}
Both the following two are well-known (for instance, see \cite{arnoldetall}). 
\begin{enumerate}
\item 
Any stable map-germ 
$(\mathbb{R}^2,0)\to (\mathbb{R}^3,0)$ is 
either a Whitney umbrella or non-singular.   
\item  
Any Legendrian stable singularity $(\mathbb{R}^2,0)\to (\mathbb{R}^3,0)$ is 
either a swallowtail or a cuspidal edge.    
\end{enumerate}     
Therefore, Theorems~\ref{theorem 1} and \ref{theorem 2} may be regarded as the fundamental theorem of calculus for stable 
map-germs $(\mathbb{R}^2,0)\to (\mathbb{R}^3,0)$ and Legendrian stable singularities 
$(\mathbb{R}^2,0)\to (\mathbb{R}^3,0)$.     
\par 
Theorems~\ref{theorem 1} and \ref{theorem 2} yield the following conjecture naturally.   
\begin{conjecture}
\begin{enumerate}
\item Let $\varphi_1, \varphi_2 : (\mathbb{R}^2,0)\to (\mathbb{R}^3,0)$ 
be two $C^\infty$ map-germs of pedal unfolding type.     
Suppose that $\varphi_1$ is $\mathcal{A}$-equivalent to $\varphi_2$.   
Then,  $\mathcal{I}(\varphi_1)$ is $\mathcal{A}$-equivalent to $\mathcal{I}(\varphi_2)$.  
 \item Let $\Phi_1, \Phi_2 : (\mathbb{R}^2,0)\to (\mathbb{R}^3,0)$ 
be two normalized Legendrian map-germs.     
Suppose that $\Phi_1$ is $\mathcal{A}$-equivalent to $\Phi_2$.   
Then,  $\mathcal{D}(\Phi_1)$ is $\mathcal{A}$-equivalent to $\mathcal{D}(\Phi_2)$.  
\end{enumerate}
\end{conjecture}
\medskip 
In \S\ref{section 2}, several preparations for the proofs of Theorems~\ref{theorem 1} 
and \ref{theorem 2} and 
the proof of Proposition~\ref{proposition 1} are given.       
Theorems~\ref{theorem 1} and \ref{theorem 2} are proved 
in \S\ref{section 3} and \S\ref{section 4} 
respectively.    
\section{Preliminaries}\label{section 2}
\subsection{Function-germs with two variables and map-germs with two variables
}
Let $\mathcal{E}_2$ be the set of $C^\infty$ 
function-germs $(\mathbb{R}^2,0)\to \mathbb{R}$ 
and 
let $m_2$ be the subset of 
$\mathcal{E}_2$ consisting of $C^\infty$ function-germs $(\mathbb{R}^2,0)\to (\mathbb{R},0)$.        
The sets $\mathcal{E}_2$ have natural $\mathbb{R}$-algebra structures.  
For a $C^\infty$ map-germ $\varphi: (\mathbb{R}^2,0)\to (\mathbb{R}^2,0)$, 
let $\varphi^*: \mathcal{E}_2\to \mathcal{E}_2$ be the 
$\mathbb{R}$-algebra homomorphism defined by $\varphi^*(u)=u\circ \varphi$.   Put 
$Q(\varphi)=\mathcal{E}_2/\varphi^*m_2\mathcal{E}_2$.        
Then, $Q(\varphi)$ is an $\mathbb{R}$-algebra.     
The following Proposition~\ref{proposition 2} is a special case of theorem (2.1) of \cite{mather4}.      
\begin{proposition} \label{proposition 2}
Let $p: (\mathbb{R}^2,0)\to (\mathbb{R},0)$ be a $C^\infty$ function-germ.     
Then, the following hold:   
\begin{enumerate}
\item The $\mathbb{R}$-algebra $Q(p(x,y), y)$ is isomorphic to $Q(x^2, y)$ if and only if 
$\frac{\partial p}{\partial x}(0,0)=0$ and $\frac{\partial^2p}{\partial x^2}(0,0)\ne 0$.    
\item The $\mathbb{R}$-algebra $Q(p(x,y), y)$ is isomorphic to $Q(x, y)$ if and only if 
$(x,y)\mapsto (p(x,y), y)$ is a germ of $C^\infty$ diffeomorphism.         
\end{enumerate}
\end{proposition}
\begin{definition}[\cite{mond}] \label{definition 2}
{\rm 
Let {$T: \mathbb{R}^2\to \mathbb{R}^2$} be the linear transformation of the form 
{$T(s,\lambda)=(-s,\lambda)$}.    
Two $C^\infty$ function-germs $p_1, p_2: (\mathbb{R}^2,0)\to (\mathbb{R},0)$ are said to be 
{$\mathcal{K}^T$-{\it equivalent}} if there exist a germ of $C^\infty$ diffeomorphism 
$h : (\mathbb{R}^2,0)\to (\mathbb{R}^2,0)$ of the form $h\circ T=T\circ h$ and 
a $C^\infty$ function-germ   
$M : (\mathbb{R}^2,(0,0))\to \mathbb{R}-\{0\}$ of the form $M\circ T=M$ 
such that $p_1\circ h(x,y)=M(x,y)p_2(x,y)$.  
}
\end{definition}
\begin{theorem}[\cite{mond}]   \label{theorem 3}
Two $C^\infty$ map-germs 
$\varphi_i: (\mathbb{R}^2,0)\to (\mathbb{R}^3,0)$ $(i=1,2)$  
of the following form 
$$
\varphi_i(x,y)=(xp_i(x^2,y),x^2,y)\quad 
\mbox{where }p_i(x^2,y)\not\in m_2^{\infty}, \quad (i=1,\; 2) 
$$
are $\mathcal{A}$-equivalent if and only if the function-germs $p_i(x^2,y)$ are 
$\mathcal{K}^T$-equivalent.   
\end{theorem}
Here, $m_2^\infty=\{q:(\mathbb{R}^2,0)\to (\mathbb{R},0)\; |\; \frac{\partial^{i+j}q }{\partial x^i\partial y^j}(0,0)=0\; 
(\forall i,j\in \{0\}\cup \mathbb{N})\}$.    
By Theorem~\ref{theorem 3} and the Malgrange preparation theorem (for instance, see \cite{arnoldetall}), 
the following holds:   
\begin{corollary} \label{corollary 1}
Two $C^\infty$ map-germs 
$\varphi_i: (\mathbb{R}^2,0)\to (\mathbb{R}^3,0)$ $(i=1,2)$  
of the following form 
$$
\varphi_i(x,y)=(n_i(x,y)p_i(x^2,y),x^2,y), 
$$
where $p_i(x^2,y)\not\in m_2^{\infty}$ 
and $n_i(x,y)$ satisfies $\frac{\partial n_i}{\partial x}(0,0)\ne 0$ for each $i\in \{1,2\}$,  
are $\mathcal{A}$-equivalent if and only if the function-germs $p_i(x^2,y)$ are 
$\mathcal{K}^T$-equivalent.  
\end{corollary}
\subsection{Map-germs of pedal unfolding type} 
Let ${\varphi}: I\times J\to \mathbb{R}^3$ be a representative 
of a given $C^\infty$ map-germ of pedal unfolding type, where 
$I, J$ be a sufficiently small intervals containing the origin of $\mathbb{R}$.     
Then, we may put $\varphi(x,y)=\left(n(x,y)p(x,y), p(x,y)\right)$.    
Put  
\[
{\Phi}(x,y)
=({\Phi}_1(x, y), {\Phi}_2(x, y),y)    
=  \left(\int_0^x n(x,y)p(x,y)dx, \int_0^x p(x,y)dx, y\right)
\]
and 
$$
\widetilde{\mu}_\Phi(x,y)={\frac{\partial }{\partial X}-n(x,y)\frac{\partial}{\partial Y}}.   
$$
Since $\widetilde{\mu}_\Phi(x,y)\ne 0$ for any $x\in I$ and $y\in J$, 
for any fixed $y\in J$ we may define the map-germ 
$L_{\Phi, y}: (\mathbb{R},0)\to T_1\mathbb{R}^2$ 
as 
$$
L_{\Phi, y}(x)=\left( \left(
\Phi_1(x, y), \Phi_2(x, y) \right) 
, \frac{\widetilde{\mu}_\Phi(x,y)}{||\widetilde{\mu}_\Phi(x,y)||}\right),         
$$
where $T_1\mathbb{R}^2$ is the unit tangent bundle of $\mathbb{R}^2$.      
Then, since $\varphi$ is a representative of a map-germ of pedal unfolding type, 
we have the following:    
\begin{lemma} \label{lemma 2.1}
For any $y\in J$,  
$L_{\Phi, y}: (\mathbb{R},0)\to T_1\mathbb{R}^2$ 
is a Legendrian lift of the map-germ $x\mapsto (\Phi_1(x,y), \Phi_2(x,y))$.    
\end{lemma}
By Lemma~\ref{lemma 2.1}, we have the following:   
\begin{lemma} \label{lemma 2.2}
For any $y\in J$, the map-germ $\widetilde{\Phi}_y: (\mathbb{R},0)\to (\mathbb{R}^2,0)$ defined by 
$\widetilde{\Phi}_y(x)=(\Phi_1(x,y), \Phi_2(x,y))$ is a Legendrian map-germ.   
\end{lemma}
\par 
\smallskip 
Next, put 
$$
\widetilde{\nu}_\Phi(x,y) 
 =  \widetilde{\mu}_\Phi(x,y)-\left(\frac{\partial \Phi_1}{\partial y}(x,y)
-n(x,y)\frac{\partial \Phi_2}{\partial y}(x,y)\right)
\frac{\partial}{\partial Z}.     
$$ 
Then,  we have the following:   
\begin{lemma} \label{lemma 2.3}
For any $x\in I$ and $y\in J$, 
$$
\widetilde{\nu}_\Phi(x,y)\cdot \frac{\partial \Phi}{\partial x}(x, y)=0,\quad 
\widetilde{\nu}_\Phi(x,y)\cdot \frac{\partial \Phi}{\partial y}(x,y)=0.   
$$
\end{lemma}
Since $\widetilde{\nu}_\Phi(x,y)\ne 0$ for any $x\in I$ and $y\in J$,    
we may define the map-germ $L_\Phi : (\mathbb{R}^2,0)\to T_1\mathbb{R}^3$ as   
\[
L_\Phi(x,y)=\left({\color{black}\Phi(x,y)}, \frac{\widetilde{\nu}_\Phi(x,y)}{||\widetilde{\nu}_\Phi(x,y)||}\right).
\] 
Then, 
by Lemma~\ref{lemma 2.3} we have the following:  
\begin{lemma} \label{lemma 2.4}
$L_\Phi: (\mathbb{R}^2,0)\to T_1\mathbb{R}^3$ is a Legendrian lift of $\Phi$.   
\end{lemma}
By Lemma~\ref{lemma 2.4} we have the following:    
\begin{lemma} \label{lemma 2.5}
$\Phi: (\mathbb{R}^2,0)\to (\mathbb{R}^3,0)$ is a Legendrian map-germ.   
\end{lemma}
%
\subsection{Normalized Legendrian map-germs}
Let $\Phi: U\to \mathbb{R}^3$ be a representative of a given normalized Legendrian map-germ 
$(\mathbb{R}^2,0)\to (\mathbb{R}^3,0)$, where $U$ is a sufficiently small 
neighborhood of the origin of $\mathbb{R}^2$.      
We assume that the origin of $\mathbb{R}^2$ is a singular point of $\Phi$.
By the condition (a) of the definition of normalized Legendrian map-germs, 
we may assume that $\Phi$ has the form $\Phi(x,y)=(\Phi_1(x,y), \Phi_2(x,y),y)$.        
Since $\Phi$ is a a representative of a 
Legendrian map-germ, we have the following:   
\begin{lemma} \label{lemma 2.6}
There exists a $C^\infty$ vector field $\nu_\Phi$ along 
$\Phi$, 
\[
\nu_\Phi(x,y)=n_1(x,y)\frac{\partial}{\partial X}
+n_2(x,y)\frac{\partial}{\partial Y}
+n_3(x,y)\frac{\partial}{\partial Z}, 
\] 
such that the following three hold:   
\begin{enumerate}
\item $n_1(x,y)\frac{\partial \Phi_1}{\partial x}(x,y)+n_2(x,y)\frac{\partial \Phi_2}{\partial x}(x,y)=0$.   
\item $n_1(x,y)\frac{\partial \Phi_1}{\partial y}(x,y)+n_2(x,y)\frac{\partial \Phi_2}{\partial y}(x,y)
+n_3(x,y)=0$. 
\item The map $L_\Phi : U\to T_1\mathbb{R}^3$ defined by 
$L_\Phi(x,y)=\left(\Phi(x,y), {\nu}_\Phi(x,y)\right)$ is an immersion.    
\end{enumerate}
\end{lemma}
By the condition (c) of the definition of normalized Legendrian map-germs, 
we have the following:   
\begin{lemma} \label{lemma 2.7}
For the vector field ${\nu}_\Phi$, 
$n_1(0,0)\ne 0$ and 
$n_2(0,0)=n_3(0,0)=0$.   
\end{lemma}
By the assertion (i) of Lemma~\ref{lemma 2.6} and Lemma~\ref{lemma 2.7}, we have the following equality as function-germs.     
\begin{equation} \label{divisible}
\frac{\partial \Phi_1}{\partial x}(x,y)=-\frac{n_2(x,y)}{n_1(x,y)}\frac{\partial \Phi_2}{\partial x}(x,y).
\end{equation}
Then, by the condition (b) of the definition of normalized Legendrian maps and the equality \eqref{divisible}, 
the following holds:   
\begin{lemma} \label{lemma 2.8}
The map-germ $\mathcal{D}(\Phi)$ maps the origin to the origin.   
\end{lemma}
Put 
\[
n(x,y)=-\frac{n_2(x,y)}{n_1(x,y)} \mbox{ and }p(x,y)=\frac{\partial \Phi_2}{\partial x}(x,y).
\]
Then, we have clearly the following:   
\begin{lemma} \label{lemma 2.9}
Both function-germs $n$ and $p$ are of class $C^\infty$ and 
the equality $\mathcal{D}(\Phi)(x,y)=(n(x,y)p(x,y), p(x,y),y)$ holds.   
\end{lemma}   
Furthrmore, we have the following:   
\begin{lemma} \label{lemma 2.10}
The function-germ $n$ satisfies that $n(0,0)=0$ and $\frac{\partial n}{\partial x}(0,0)\ne 0$.  
\end{lemma}
\underline{\it Proof.}\qquad  
By Lemma~\ref{lemma 2.7}, we have that $n(0,0)=0$.      
Suppose that $\frac{\partial n}{\partial x}(0,0)=0$.     
Then, by differentiating both side of the equality in the assertion (ii) of Lemma~\ref{lemma 2.6} with respect to $x$, 
we have the following equality:   
\[
n_1(0,0)\frac{\partial^2\Phi_1}{\partial x\partial y}(0,0)+\frac{\partial n_3}{\partial x}(0,0)=0.   
\] 
Since we have assumed $\frac{\partial n}{\partial x}(0,0)=0$, 
we have that $\frac{\partial n_2}{\partial x}(0,0)=0$.     
Thus and since $\Phi$ is a normalized Legendrian map-germ 
such that the origin of $\mathbb{R}^2$ is a singular point of $\Phi$, 
we have that $\frac{\partial n_3}{\partial x}(0,0)\ne 0$.     
Thus, we have that 
$\frac{\partial^2\Phi_1}{\partial x\partial y}(0,0)\ne 0$.      
Hence, by the condition (b) of the definition of normalized Legendrian maps, 
Lemma~\ref{lemma 2.7} and the equality \eqref{divisible}, 
we have a contradiction.     
\hfill $\Box$ 
\par 
\smallskip 
\begin{definition} \label{definition 3}
{\rm 
Let $\Phi: (\mathbb{R}^2, 0)\to (\mathbb{R}^3,0)$ be a Legendrian map-germ and let 
$\nu_\Phi$ be a unit normal vector field of $\Phi$ given in the definition of Legendrian map-germs.      
The $C^\infty$ function-germ $LJ_{\Phi}: (\mathbb{R}^2,0)\to \mathbb{R}$ defined by the following is 
called the {\it Legendrian-Jacobian} of $\Phi$.    
\[
LJ_{\Phi}(x,y)= 
\det 
\left(
\frac{\partial \Phi}{\partial x}(x,y), \frac{\partial \Phi}{\partial y}(x,y), {\nu}_\Phi(x,y)
\right).    
\]
}
\end{definition}
Note that if $\nu_\Phi$ satisfies the conditions of unit normal vector field of $\Phi$, 
then $-\nu_\Phi$ also satisfies 
them.   
Thus, the sign of $LJ_{\Phi}(x,y)$ depends on the particular choice of unit normal vector field $\nu_\Phi$.     
The Legendrian Jacobian of $\Phi$ is called also 
the {\it signed area density function} 
(for instance, see \cite{sajiumeharayamada2}).       
Although it seems reasonable to call $LJ_{\Phi}$ the area density function from the viewpoint of 
investigating the singular surface $\Phi(U)$ ($U$ 
is a sufficiently small neighborhood of the origin of $\mathbb{R}^2$), 
it seems reasonable to call it the Legendrian Jacobian from the viewpoint of 
investigating the singular map-germ $\Phi$.      
Hence, we call $LJ_{\Phi}$ the Legendrian Jacobian of $\Phi$ in this paper.    
\par 
\medskip 
Let $\Phi: (\mathbb{R}^2, 0)\to (\mathbb{R}^3,0)$ be a normalized Legendrian map-germ  
and $\nu_\Phi$ is a unit normal vector field of $\Phi$.       
\begin{eqnarray*}
\Phi(x,y) & = & (\Phi_1(x,y), \Phi_2(x,y), y), \\
\nu_\Phi(x,y) & = & n_1(x,y)\frac{\partial}{\partial X}+n_2(x,y)\frac{\partial}{\partial Y}+n_3(x,y)\frac{\partial}{\partial Z}.    
\end{eqnarray*}
By Lemma~\ref{lemma 2.7}, we may put 
\[
\widetilde{\nu}_\Phi(x,y)=\frac{\partial}{\partial X}+\frac{n_2(x,y)}{n_1(x,y)}\frac{\partial}{\partial Y}
+\frac{n_3(x,y)}{n_1(x,y)}\frac{\partial}{\partial Z}.    
\]
\begin{lemma} \label{lemma 2.11}
The Legendrian Jacobian $LJ_{\Phi}$ is expressed as follows:
\[
LJ_\Phi(x,y)=\frac{\frac{\partial \Phi_2}{\partial x}(x,y)}{n_1(x,y)}.   
\]
\end{lemma}
\underline{\it Proof.}\qquad     
Calculations show that 
\[
 \frac{\partial \Phi}{\partial x}(x,y)\times\frac{\partial \Phi}{\partial y}(x,y) 
  =  \frac{\partial \Phi_2}{\partial x}(x,y)\widetilde{\nu}_\Phi(x,y) 
\]
where the cross in the center stands for the vector product.      
It follows $LJ_\Phi(x,y)=\frac{\partial \Phi_2}{\partial x}(x,y)/n_1(x,y)$.    
\hfill $\Box$
\subsection{Proof of Proposition~\ref{proposition 1}} 
{ } \quad 
\\ 
\noindent 
\underline{\it Proof of the assertion (i) of Proposition~\ref{proposition 1}.}\qquad     
\par 
Put $\mathcal{I}(\varphi)=\Phi(x,y)=\left(\Phi_1(x,y), \Phi_2(x,y), y\right)$.      
Then, by Lemma~\ref{lemma 2.5}, $\Phi$ is a Legendrian map-germ.    
Thus, in order to complete the proof of the assertion (i) of Proposition~\ref{proposition 1}, it is sufficient to 
show the conditions (b), (c) of the definition of normalized Legendrian map-germs are satisfied.      
\par 
Put $\varphi(x,y)=(n(x,y)p(x,y), p(x,y),y)$.   Then, by the definition of map-germs of pedal unfolding type, 
we have that $n(0,0)=0$ and $p(0,0)=0$.     It follows that $\frac{\partial \Phi_2}{\partial x}(0,0)=p(0,0)=0$.     
Thus, the condition (b) is satisfied.       
By Lemma~\ref{lemma 2.4},  the following $L_\Phi$ is a germ of Legendrian lift of $\Phi$.      
\[
L_\Phi(x,y)=\left(\Phi(x,y), \frac{\widetilde{\nu}_\Phi(x,y)}{||\widetilde{\nu}_\Phi(x,y) ||}\right),  
\]
where 
$\widetilde{\nu}_\Phi(x,y)=\frac{\partial}{\partial X}-n(x,y)\frac{\partial}{\partial Y} 
- \left(\frac{\partial \Phi_1}{\partial y}(x,y)-n(x,y)\frac{\partial \Phi_2}{\partial y}(x,y)\right)\frac{\partial}{\partial Z}$.     
Since $n(0,0)=0$ and $\frac{\partial \Phi_1}{\partial y}(0,0)=\int_0^0\frac{\partial np}{\partial y}(x,0)dx=0$, we have 
\[
\frac{\widetilde{\nu}_\Phi(0,0)}{||\widetilde{\nu}_\Phi(0,0) ||}=\frac{\partial}{\partial X}.   
\]
 Thus the condition (c) is satisfied.   
 \hfill $\Box$
\par 
\smallskip 
\noindent  
\underline{\it Proof of the assertion (ii) of Proposition~\ref{proposition 1}.}\qquad     
\par 
The assertion (ii) of Proposition~\ref{proposition 1} 
follows from Lemmas~\ref{lemma 2.8}, \ref{lemma 2.9} and \ref{lemma 2.10}.   
\hfill $\Box$ 
\section{Proof of Theorem~\ref{theorem 1}}  \label{section 3}   
Suppose that both $\mathcal{I}: \mathcal{W}\to \mathcal{S}$ and 
$\mathcal{D}: \mathcal{S}\to \mathcal{W}$ are well-defined.    
Then, by the fundamental theorem of calculus, the following hold:   
\begin{eqnarray*}
\mathcal{D}\circ \mathcal{I}(\varphi) & = & \varphi\quad \mbox{ for any }\varphi\in \mathcal{W},  \\ 
\mathcal{I}\circ \mathcal{D}(\Phi) & = & \Phi\quad \mbox{ for any }\Phi\in \mathcal{S}. \\ 
\end{eqnarray*}  
Thus, both $\mathcal{I}$ and $\mathcal{D}$ are bijective.    
Therefore, in order to complete the proof, it is sufficient to show that 
both $\mathcal{I}$ and $\mathcal{D}$ are well-defined.  
\par 
\medskip
\noindent 
\underline{\it Proof that $\mathcal{I} : \mathcal{W}\to \mathcal{S}$ is well-defined.}\qquad     
\par 
Let $\varphi(x,y)=(n(x,y)p_\varphi(x,y), p_\varphi(x,y),y)$ be an element  of $\mathcal{W}$.        
Put $\Phi=\mathcal{I}(\varphi)$.     
Then, 
$\Phi$ is a normalized Legendrian map-germ by Proposition~\ref{proposition 1} in \S\ref{intro}.        
Let $g$ be the Whitney umbrella of pedal unfolding type \eqref{g} defined in \S\ref{intro}:  
$$g(x,y)=(xp_g(x,y), p_g(x,y), y)=(x(x^2+y), -x^2-y,y).$$      
\begin{lemma} \label{lemma 3.1}
There exists a germ of $C^\infty$ diffeomorphism 
$h: (\mathbb{R}^2,0)\to (\mathbb{R}^2,0)$ such that 
$h$ has the form $h(x,y)=(h_1(x,y), h_2(y))$ and $p_\varphi\circ h(x,y)$ is $x^2+y$ or $-(x^2+y)$.     
\end{lemma}       
\underline{\it Proof.}\qquad     
Since $\varphi$ is 
a Whitney umbrella of pedal unfolding type, 
we have the following:   
\[ 
Q(p_\varphi(x,y), y)
\cong 
Q(\varphi)\cong Q(g)
\cong 
Q(x^2, y).   
\]    
Thus, we may put $p_\varphi(x,0)=a_2x^2+o(x^2)$ $(a_2\ne 0)$ by Proposition~\ref{proposition 2} 
in \S\ref{section 2}.    
By the Morse lemma with parameters (for instance, see \cite{brucegiblin}), 
there exists a germ of $C^\infty$ diffeomorphism 
$h: (\mathbb{R}^2,0)\to (\mathbb{R}^2,0)$ such that 
$h$ has the form $h(x,y)=(h_{1}(x,y), h_{2}(y))$ and 
$p_\varphi\circ h(x,y)=\pm (x^2+ q(y))$ 
by a certain $C^\infty$ function-germ $q:(\mathbb{R},0)\to (\mathbb{R},0)$.  
Since $\varphi$ is $\mathcal{A}$-equivalent to $g$, 
by Corollary~\ref{corollary 1} in \S\ref{section 2}, $\pm (x^2+q(y))$ is $\mathcal{K}^T$-equivalent to $p_g$ and thus 
$q: (\mathbb{R},0)\to (\mathbb{R},0)$ is a germ of $C^\infty$ diffeomorphism.     
Hence, Lemma~\ref{lemma 3.1} follows.       
\hfill $\Box$ 
\par 
\medskip 

\begin{lemma} \label{lemma 3.2}
The normalized Legendrian map-germ $\Phi$ 
is a swallowtail.    
\end{lemma}
\underline{\it Proof.}\qquad     
Put $G=\mathcal{I}(g)$.     
Then, $G$ has the form \eqref{G} in \S 1 which is a normalized swallowtail.   
Since $G$ is normalized, 
$\frac{\partial}{\partial x}$ 
is the null vector field for $G$ defined in \cite{kokuburossmansajiumeharayamada, sajiumeharayamada}, 
that is, $\frac{\partial G}{\partial x}(x,y)=0$ holds for any $(x,y)$ which is a singular point of $G$.          
Thus and since $G$ is a swallowtail, the following three hold by corollary 
2.5 of \cite{sajiumeharayamada}.    
\begin{enumerate}
\item $LJ_G(0,0)=\frac{\partial LJ_G}{\partial x}(0,0)=0$, 
\item $\frac{\partial^2 LJ_G}{\partial x^2}(0,0)\ne 0$.   
\item $Q(LJ_G, \frac{\partial LJ_G}{\partial x})\cong Q(x,y)$.    
\end{enumerate}   
\par 
On the other hand, by Lemmas~\ref{lemma 2.11} and \ref{lemma 3.1}, 
there exist a germ of $C^\infty$ diffeomorphism 
$h: (\mathbb{R}^2,0)\to (\mathbb{R}^2,0)$ and 
a $C^\infty$ function-germ $\xi: (\mathbb{R}^2,0)\to \mathbb{R}$ 
such that $h$ has the form $h(x,y)=(h_1(x,y), h_2(y))$, $\xi(0,0)\ne 0$ 
and the following hold:   
\[
LJ_\Phi\circ h(x,y)=
\xi(x,y)LJ_G(x,y).     
\]
Hence and since $\frac{\partial}{\partial x}$ 
is the null vector field for $\Phi$ (this is because $\Phi$ is normalized), 
the following three hold for $LJ_\Phi$.   
\begin{enumerate}
\item $LJ_\Phi(0,0)=\frac{\partial LJ_\Phi}{\partial x}(0,0)=0$, 
\item $\frac{\partial^2 LJ_\Phi}{\partial x^2}(0,0)\ne 0$, 
\item $Q(LJ_\Phi, \frac{\partial LJ_\Phi}{\partial x})\cong Q(x,y)$.    
\end{enumerate}   
Hence, $\Phi$ is a swallowtail 
by corollary 2.5 of \cite{sajiumeharayamada}.     
\hfill $\Box$    
\par 
\medskip
\noindent 
\underline{\it Proof that $\mathcal{D} : \mathcal{S}\to \mathcal{W}$ is well-defined.}\qquad     
\par 
Let $\Phi$ be an element of $\mathcal{S}$.       
Then, by Proposition~\ref{proposition 2} in \S\ref{section 2}, $\mathcal{D}(\Phi)$ is of pedal unfolding type.    
\begin{lemma}     \label{lemma 3.3} 
For the Legendrian Jacobian $LJ_\Phi$, the following three hold:
\begin{enumerate}
\item $LJ_\Phi(0,0)=\frac{\partial LJ_\Phi}{\partial x}(0,0)=0$, 
\item $\frac{\partial^2 LJ_\Phi}{\partial x^2}(0,0)\ne 0$.   
\item $Q(LJ_\Phi, \frac{\partial LJ_\Phi}{\partial x})\cong Q(x,y)$.    
\end{enumerate}
\end{lemma}   
\underline{\it Proof.}\qquad  
Since $\Phi$ is normalized, 
$\frac{\partial}{\partial x}$ is the null vector field.    
Thus and since $\Phi$ is a swallowtail, 
Lemma~\ref{lemma 3.3} follows from corollary 2.5 of \cite{sajiumeharayamada}.   
\begin{lemma} \label{lemma 3.4}
For the $\Phi$, the map-germ of pedal unfolding type 
$\mathcal{D}(\Phi)$ is a Whitney umbrella.   
\end{lemma}
\underline{\it Proof.}\qquad  
Since $\mathcal{D}(\Phi)$ is of pedal unfolding type, 
there exists a $C^\infty$ function-germ 
$n: (\mathbb{R}^2,0)\to (\mathbb{R},0)$ such that 
$\frac{\partial n}{\partial x}(0,0)\ne 0$ and 
$\frac{\partial \Phi_1}{\partial x}(x,y)=n(x,y)\frac{\partial \Phi_2}{\partial x}(x,y)$ 
where $\Phi(x,y)=(\Phi_1(x,y),$ $\Phi_2(x,y),y)$.      
Put $p_\varphi=\frac{\partial \Phi_2}{\partial x}$.     
Then, by Lemmas~\ref{lemma 2.11}, \ref{lemma 3.3} and the Morse lemma with parameters, 
there exists a germ of $C^\infty$ diffeomorphism 
$h: (\mathbb{R}^2,0)\to (\mathbb{R}^2,0)$ such that 
$h$ has the form $h(x,y)=(h_1(x,y), h_2(y))$ and $p_\varphi\circ h(x,y)$ is $x^2+y$ or $-(x^2+y)$.     
Therefore, by Corollary~\ref{corollary 1} in \S\ref{section 2}, $\mathcal{D}(\Phi)$ is $\mathcal{A}$-equivalent to $g$.      
\hfill $\Box$ 
\par 
\smallskip 
\section{Proof of Theorem~\ref{theorem 2}}     \label{section 4}
As same as Theorem~\ref{theorem 1}, it is sufficient to show that 
both $\mathcal{I}: \mathcal{N}\to \mathcal{C}$ and 
$\mathcal{D}: \mathcal{C}\to \mathcal{N}$ are well-defined.  
\par 
\medskip
\noindent 
\underline{\it Proof that $\mathcal{I} : \mathcal{N}\to \mathcal{C}$ is well-defined.}\qquad     
\par 
Let $\varphi(x,y)=(n(x,y)p_\varphi(x,y), p_\varphi(x,y),y)$ be an element  of $\mathcal{N}$.     
Put $\Phi=\mathcal{I}(\varphi)$.     
Then, since $\varphi$ is of pedal unfolding type, $\Phi$ is a normalized 
Legendrian map-germ by Proposition~\ref{proposition 1} in \S\ref{intro}.       
Let $g$ be the non-singular map-germ of pedal unfolding type defined by 
$
g(x,y)= (x^2, x, y)$.     
\begin{lemma} \label{lemma 4.1}
There exists a germ of $C^\infty$ diffeomorphism 
$h: (\mathbb{R}^2,0)\to (\mathbb{R}^2,0)$ such that 
$h$ has the form $h(x,y)=(h_1(x,y), h_2(y))$ and $p_\varphi\circ h(x,y)=x$ holds.     
\end{lemma}       
\underline{\it Proof.}\qquad     
Since $\varphi$ is non-singular and of pedal unfolding type, we have the following:   
\[ 
Q(p_\varphi(x,y), y)
\cong 
Q(\varphi)\cong Q(g)
\cong 
Q(x,y).
\]    
Thus, $(p_\varphi(x,y), y)$ is a germ of $C^\infty$ diffeomorphism 
by Proposition~\ref{proposition 2} in \S\ref{section 2}.         
From the form of $(p_\varphi(x,y), y)$, its inverse map-germ 
$h: (\mathbb{R}^2,0)\to (\mathbb{R}^2,0)$ has the form $h(x,y)=(h_{1}(x,y), h_{2}(y))$.   
Since $h$ is the inverse map-germ of $(p_\varphi(x,y), y)$, 
it follows that $p_\varphi\circ h(x,y)=x$.    
\hfill $\Box$ 
\begin{lemma} \label{lemma 4.2}
The normalized Legendrian map-germ $\Phi$ is a cuspidal edge.    
\end{lemma}
\underline{\it Proof.}\qquad     
Since $\Phi$ is normalized, $\frac{\partial}{\partial x}$ is the null vector field for $\Phi$.      
By Lemmas~\ref{lemma 2.11}, \ref{lemma 4.1}, we have that $\frac{\partial LJ_\Phi}{\partial x}(0,0)\ne 0$.     Thus,  
thus the null vector field $\frac{\partial}{\partial x}$ is transverse to $\{(x,y)\; |\; LJ_\Phi(x,y)=0\}$ 
at $(0,0)\in \mathbb{R}^2$.         
Hence, $\Phi$ is a cuspidal edge 
by proposition 1.3 of \cite{kokuburossmansajiumeharayamada}.     
\hfill $\Box$ 
\par 
\medskip
\noindent 
\underline{\it Proof that $\mathcal{D} : \mathcal{C}\to \mathcal{N}$ is well-defined.}\qquad     
\par 
Let $\Phi$ be an element of $\mathcal{C}$.      Then, by Proposition~\ref{proposition 1} 
in \S\ref{intro}, 
$\mathcal{D}(\Phi)$ is of pedal unfolding type.    
\begin{lemma}    \label{lemma 4.3}  
For the Legendrian Jacobian $LJ_\Phi$, 
two properties $LJ_\Phi(0,0)=0$ and 
$\frac{\partial LJ_\Phi}{\partial x}(0,0)\ne 0$ hold.    
\end{lemma}   
\underline{\it Proof.}\qquad  
Since $\frac{\partial}{\partial x}$ is the null vector field for $\Phi$ 
and $\Phi$ is a cuspidal edge, 
Lemma~\ref{lemma 4.3} follows from corollary 2.5 of \cite{sajiumeharayamada}.   
\begin{lemma} \label{lemma 4.4}
For the $\Phi$, the map-germ of pedal unfolding type 
$\mathcal{D}(\Phi)$ is non-singular and 
$\mathcal{D}(\Phi)(0,0)=(0,0,0)$.   
\end{lemma}
\underline{\it Proof.}\qquad  
Since $\mathcal{D}(\Phi)$ is of pedal unfolding type, 
there exists a $C^\infty$ function-germ 
$n: (\mathbb{R}^2,0)\to (\mathbb{R},0)$ such that 
$\frac{\partial n}{\partial x}(0,0)\ne 0$ and 
$\frac{\partial \Phi_1}{\partial x}(x,y)=n(x,y)\frac{\partial \Phi_2}{\partial x}(x,y)$ 
where $\Phi(x,y)=(\Phi_1(x,y),$ $\Phi_2(x,y),y)$.       
Put $p_\varphi=\frac{\partial \Phi_2}{\partial x}$.     
Then, by Lemmas~\ref{lemma 2.11} and \ref{lemma 4.3}, 
the map-germ $(x,y)\mapsto (p_\varphi(x,y), y)$ is a germ of $C^\infty$ diffeomorphism.   
Thus, $\mathcal{D}(\Phi)$ is non-singular.
\hfill $\Box$ 
\par 
\smallskip 
\subsection*{Acknowledgement}
The author 
would like to thank Kentaro Saji for valuable comments.   

{\color{black}
\begin{flushright}

 

\end{flushright}
}
\end{document}